\documentclass[12pt,oneside]{article}
\usepackage[cp866]{inputenc}   
\usepackage{color}
\usepackage{amsmath}
\usepackage{amssymb}
\usepackage{amsthm}
\usepackage[english]{babel}

\newtheorem{theo}{\bf Theorem}

\title{{\bf {\bf \normalsize Alexander L. Gavrilyuk\footnote{
~e-mail:~\texttt{alexander.gavriliouk@gmail.com}\newline
N.N. Krasovsky Institute of Mathematics and Mechanics of Russian Academy of Sciences, Ekaterinburg, Russia.},~
Ivan Y. Mogilnykh\footnote{
~e-mail:~\texttt{ivmog84@gmail.com}\newline
S.L. Sobolev Institute of Mathematics of Russian Academy of Sciences, Novosibirsk, Russia.}\\}
\normalsize \large ON THE GODSIL -- HIGMAN NECESSARY CONDITION FOR EQUITABLE 
PARTITIONS OF ASSOCIATION SCHEMES}}
\begin{document}

\maketitle
\section{Introduction}
\hspace{5mm}

Equitable partitions of association schemes are often related to
important and interesting combinatorial and geometric objects such as
combinatorial designs, orthogonal arrays and the Cameron --
Liebler line classes in projective geometries. Thereby the question of existence
of these structures in association schemes is very difficult in general.

In this short note, we show that a necessary condition for
equitable partitions recently proposed by C. Godsil in his
monograph \cite{God10} is not stronger than the well-known Lloyd
theorem (see Theorem 1 in Section 2). In the next section, we
recall some basic definitions and notions. Section 3 contains the
proof of our result (Theorem 2).

\section{Association schemes and their equitable partitions}

Let $V$ be a finite set of size $v$ and ${\mathbf C}^{V\times V}$
be the set of matrices over ${\mathbf C}$ with rows and columns
indexed by $V$. Let ${\mathcal R}=\{R_0,R_1,\ldots,R_{d}\}$ be a
set of non-empty subsets of $V\times V$. For $i=0,\ldots,d$, let
$A_i\in {\mathbf C}^{V\times V}$ be the adjacency matrix of the
graph $(V,R_i)$. The pair ${\mathcal A}=(V,{\mathcal R})$ is said to be
an {\it association scheme} with $d$ classes and vertex set $V$
if the following properties hold:
\begin{enumerate}
\item[(1)] $A_0=I$, the identity matrix,
\item[(2)] $\sum_{i=0}^d A_i=J$, where $J$ is the all ones matrix,
\item[(3)] $A_i^{\rm T}\in \{A_0,\ldots,A_d\}$, for every $i=0,\ldots,d$,
\item[(4)] $A_iA_j$ is a linear combination of $A_0,\ldots,A_d$,
for all $i,j=0,\ldots,d$.
\end{enumerate}

The matrix algebra ${\mathbf C}[{\mathcal A}]$ over ${\mathbf C}$
generated by $A_0,\ldots,A_d$ is called the {\it Bose -- Mesner algebra}
of ${\mathcal A}$. It now follows from properties (1)-(4) that
${\mathbf C}[{\mathcal A}]$ has a basis consisting of the
matrices $A_0,\ldots,A_d$ and its dimension is $d+1$. We say that
${\mathcal A}$ is commutative if ${\mathbf C}[{\mathcal A}]$ is
commutative, and that ${\mathcal A}$ is symmetric if the matrices
$A_i$ are symmetric. Clearly, a symmetric association scheme is
commutative.

In the paper, we assume that ${\mathcal A}$ is a symmetric
association scheme with $d$ classes.

Since ${\mathcal A}$ is commutative, it follows that the matrices
$A_0$, $A_1$, $\ldots$, $A_d$ are simultaneously diagonalized by
an appropriate unitary matrix. This means that ${\mathbf
C}^{V}$ is decomposed as an orthogonal direct sum of $d+1$
maximal common eigenspaces of $A_0$, $A_1$, $\ldots$,
$A_d$:
$${\mathbf C}^{V}=W_0\oplus W_1\oplus \ldots W_d,$$
and, for every $0\leq j\leq d$, define $E_j\in {\mathbf C}^{V\times V}$
to be the orthogonal projection onto $W_j$.
Let $\{\overline{w}_{j\ell}\mid 1\leq\ell\leq {\rm dim}(W_j)\}$ be an orthonormal
basis of $W_j$, $0\leq j\leq d$.

Note \cite{Bra89} that the matrices $E_j$ form another basis for
${\mathcal A}$ consisting of the primitive idempotents of
${\mathcal A}$, i.e., $E_iE_j = \delta_{ij}E_i$ and
$\sum_{j=0}^dE_j=I$.

For the two basises $A_0,\ldots,A_d$ and $E_0,\ldots,E_d$ of ${\mathbf C}[{\mathcal A}]$,
the change-of-basis matrices $P$ and $Q$ are defined by
$$A_i=\sum_{j=0}^dP_{ji}E_j,~~E_i=\frac{1}{v}\sum_{j=0}^dQ_{ji}A_j,$$
where, in fact, $P_{ji}$ is an eigenvalue of $A_i$ on the eigenspace $W_j$.
It now follows that
$$PQ=vI.$$

The numbers $v_i=P_{0i}$, $0\leq i\leq d$, are called the {\it
valencies} of the scheme. The numbers $f_j={\rm tr}(E_j)={\rm
rank}(E_j)={\rm dim}(W_j)$, $0\leq j\leq d$, are called the {\it
multiplicities} of the scheme. By \cite[Lemma 2.2.1(iv)]{Bra89} 
the following relation holds:
\begin{equation}\label{crossPQ}
\frac{Q_{ij}}{f_j}=\frac{P_{ji}}{v_i}.
\end{equation}

The distance-regular graphs provide important but not only
examples of symmetric association schemes. For more results and
background on distance-regular graphs and association schemes, we
refer the reader to \cite{Ban84}, \cite{Bra89}, \cite{God04}.
\smallskip

Let $\pi$ be a partition of $V$ with $t$ cells $C_1,C_2,\ldots,C_t$.
The characteristic matrix $H$ of $\pi$ is the $v\times t$-matrix whose columns
are the characteristic vectors of ${\mathbf C}^{V}$ of the cells of $\pi$.

We say that a partition $\pi$ of ${\mathcal A}$ is {\it equitable}
\cite{God10} if the column space of $H$ is ${\mathbf C}[{\mathcal
A}]$-invariant, i.e., for every $A\in {\mathbf C}[{\mathcal A}]$,
there exists a $t\times t$-matrix $B$ usually called the {\it
quotient} matrix of $A$ such that
$$AH=HB.$$

One can give the following equivalent combinatorial definition of equitable partition:
$\pi$ is equitable if, for all $0\leq i\leq d$, $1\leq j,k\leq t$ and every vertex
$x\in C_k$, there exist exactly $n_{ij}^k$ vertices $y\in C_j$ such that $(x,y)\in R_i$.
We define the $t\times t$-matrices $N_i$, $0 \leq i\leq d$, by
$$(N_i)_{k,j}=n_{ij}^k,$$
and it is easy to see that, in fact, $N_i$ is the quotient matrix of $A_i$, i.e.,
$A_iH=HN_i.$
\smallskip

Let $C$ be a subset of the vertex set of a graph $\Gamma$. The covering radius
$\rho_C$ of $C$ is defined to be $\rho_C:={\rm max}\{{\rm d}(x,C)\mid x\in \Gamma\}$,
where ${\rm d}(x,C):={\rm min}\{{\rm d}(x,y)\mid y\in C\}$, and ${\rm d}(x,y)$
is the usual graph distance. For $i=0,\ldots,\rho_C$, define $\Gamma_i(C)$ to be
the set of vertices that are at distance $i$ from $C$.
The partition $\{C = \Gamma_0(C),\Gamma_1(C),\ldots,\Gamma_{\rho_C}(C)\}$ is
referred to as the {\it distance partition} of the vertex set of $\Gamma$
with respect to $C$.

For an association scheme ${\mathcal A}=(V,{\mathcal R})$, assume
that a graph $\Gamma=(V,R_i)$ is distance-regular for some $i$. An
important type of equitable partitions is provided by completely
regular codes. A vertex subset $C$ of $\Gamma$ is called a {\it
completely regular code}  if the distance partition with respect
to $C$ is equitable. Clearly, a single vertex of the graph
$\Gamma$ is a completely regular code.

For instance, the completely regular codes of the Hamming graphs give 
rise to orthogonal arrays, of the Johnson graphs --- to combinatorial
designs \cite{Del73}, \cite{Martphd}, and of the Grassmann graphs
of diameter 2 --- to the Cameron -- Liebler line classes in the projective
geometry of dimension 3 \cite{Van11}. 

Another interesting type of equitable partitions is provided by
the partition of the vertex set of antipodal distance-regular
graphs into antipodal classes \cite{Martphd}.
\smallskip

The following result is well known in Algebraic Combinatorics and
is sometimes referred to as Lloyd's theorem \cite{God04}.

\begin{theo} If $\pi$ is an equitable partition of symmetric association scheme
${\mathcal A}$, then, for any matrix $A\in {\mathcal A}$, the characteristic
polynomial of its quotient matrix $B$ divides the characteristic polynomial of $A$.
\end{theo}

In other words, Theorem 1 states that every eigenvalue of the quotient matrix $B$
is an eigenvalue of $A$, and its multiplicity as an eigenvalue of $B$ is not greater
than its multiplicity as an eigenvalue of $A$.

\section{Main result}

Following Godsil \cite{God10}, let us define a complex inner
product on ${\mathbf C}^{V\times V}$ by
$$\langle M,N\rangle={\rm tr}(M^{*}N)={\rm sum}(\overline{M}\circ N),$$
where ${\rm sum}(M)$ denotes the sum of the entries of $M$, and
$\circ$ denotes the Schur multiplication of matrices. Then the
basis $A_0,A_1,\ldots,A_d$ of the Bose -- Mesner algebra ${\mathbf
C}[{\mathcal A}]$ is orthogonal with respect to the inner product.

Further, for a matrix $M\in {\mathbf C}^{V\times V}$, let
$\hat{M}$ denote its orthogonal projection onto ${\mathbf C}[{\mathcal A}]$, i.e.,
$$\hat{M}=\sum_{i=0}^d \frac{\langle A_i,M\rangle}{\langle A_i,A_i\rangle}A_i=
\sum_{j=0}^d \frac{\langle E_j,M\rangle}{\langle E_j,E_j\rangle}E_j.$$

It follows from \cite[Theorem 2.2.1]{God10} that
\begin{equation}\label{SeidelEq}
\hat{M}=\sum_{i=0}^d \frac{\langle A_i,M\rangle}{vv_i}A_i=
\sum_{j=0}^d \frac{\langle E_j,M\rangle}{m_j}E_j.
\end {equation}

Let $P$ be a $v\times v$ permutation matrix. Then $P$ is an
automorphism of ${\mathcal A}$ if it commutes with each matrix of
${\mathbf C}[{\mathcal A}]$. G. Higman derived the following
necessary condition for $P$ to be an automorphism of ${\mathcal
A}$, see also \cite{Cam}.

Let $P$ be an automorphism of ${\mathcal A}$ and $\sigma$ denote
the corresponding  permutation associated with $P$. Define
$\alpha_i(\sigma)$ to be the number of vertices $x$ of the vertex
set of ${\mathcal A}$ such that $(x,\sigma(x))\in R_i$. Note that
$\alpha_i(\sigma)=\langle A_i,P\rangle$. Then using equality
(\ref{SeidelEq}):
$$\hat{P}=\sum_{i=0}^d \frac{\alpha_i(\sigma)}{vv_i}A_i=
\sum_{i=0}^d \frac{\langle P,E_i\rangle}{f_i}E_i,$$ and, further,
exploiting this equation, one can show (see \cite{God10}) that a
number
$$\langle P,E_j\rangle = \frac{f_j}{v}\sum_{i=0}^d\frac{P_{ji}}{v_i}\alpha_i(\sigma)$$
must be an algebraic integer.

This condition is widely used in a study of feasible automorphisms
of distance-regular graphs \cite{Makh10}.

Now let $F$ be a projection matrix (i.e., $F^2=F$) that commutes 
with ${\mathbf C}[{\mathcal A}]$. In his monograph \cite{God10}, 
C. Godsil notices that the arguments that prove the Higman condition 
also yield that
\begin{equation}\label{conditionGodsil}
\langle F,E_j\rangle = \frac{f_j}{v}\sum_{i=0}^d\frac{P_{ji}}{v_i}\langle F,A_i\rangle
\end{equation}
must be a non-negative integer and this observation ``could be used
to show that certain equitable partitions do not exist''.

Let us define a projection matrix related to an equitable
partition. Let $\pi$ be a partition of ${\mathcal A}$ with
characteristic matrix $H$ and the cells $C_1,\ldots,C_t$. Then
$D=H^{\rm T}H$ is the diagonal matrix such that $D_{i,i}=|C_i|$.
Now the matrix $F=HD^{-1}(HD^{-1})^{\rm T}$ represents an
orthogonal projection onto the column space of $H$. For an
appropriate ordering of $X$, one can see that $F$ has the
following block form:
$$ F = \left(
\begin{array}{cccc}
  \frac{1}{|C_1|}J_{|C_1|} & O_{|C_1|\times |C_2|} & \ldots & O_{|C_1|\times |C_t|} \\
  O_{|C_2|\times |C_1|} & \frac{1}{|C_2|}J_{|C_2|} & \ldots & O_{|C_2|\times |C_t|} \\
  \ldots & \ldots & \ldots & \ldots \\
  O_{|C_t|\times |C_1|} & \ldots & O_{|C_t|\times |C_{t-1}|} & \frac{1}{|C_t|}J_{|C_t|}\\
\end{array}\right),$$

\noindent where $J$ is the all ones square matrix of proper order.
Note that $\pi$ is equitable if and only if $F$ commutes with
${\mathbf C}[{\mathcal A}]$.

Now let $\pi$ be a putative equitable
partition of ${\mathcal A}$ with quotient matrices
$N_0,\ldots,N_d$, characteristic matrix $H$, and $F$ be the
corresponding projection matrix. For a left eigenvector
$\overline{w}_{j\ell}$ of $A_i$ with eigenvalue $P_{ji}$ we have
$\overline{w}_{j\ell}A_iH=P_{ji}\overline{w}_{j\ell}H=\overline{w}_{j\ell}HN_i$
so that $\overline{w}_{j\ell}H$ is a left eigenvector of $N_i$
with the same eigenvalue $P_{ji}$ if $\overline{w}_{j\ell}H\ne
\overline{0}$. Thus, if a space $W_jH$ is non-zero then it is a subspace 
of an eigenspace of $N_i$. Since $H$ has rank $t$, it follows that there are 
exactly $t$ linearly independent vectors of type $\overline{w}_{j\ell}H$ 
for some $0\leq j\leq d$, which are left eigenvectors of $N_i$. 

Let $m_j$ denote the dimension of $W_jH$, $1\leq j\leq d$. It is now
easily seen that, for each $k=0,\ldots,t$, the spectrum of $N_k$
consists of the numbers $P_{0k},\ldots,P_{dk}$ (some of them might
be equal) with multiplicities $m_0$, $m_1$, $\ldots$, $m_d$
respectively.


\begin{theo} The following equality holds:
$$\langle F,E_j\rangle = m_j,$$ where $m_j$ is the multiplicity of the eigenvalue
$P_{ji}$ as an eigenvalue of $N_i$
 on the eigenspace $W_jH$, for
every $0 \leq i\leq d$.

\end{theo}

\noindent{\it Proof.} First of all, we note that, for $0\leq i\leq d$,
$$\langle F,A_i\rangle = {\rm tr}(FA_i)=n_{i1}^1+n_{i2}^2+\ldots+n_{it}^t={\rm tr}(N_i),$$
and therefore $\langle F,A_i\rangle$ is equal to the sum of all eigenvalues of $N_i$.

We now have
$$\langle F,E_j\rangle =
\frac{f_j}{v}\sum_{i=0}^d\frac{P_{ji}}{v_i}\langle F,A_i\rangle =
\frac{f_j}{v}\sum_{i=0}^d\frac{P_{ji}}{v_i}{\rm tr}(N_i) =
\frac{f_j}{v}\sum_{i=0}^d\frac{P_{ji}}{v_i}\sum_{k=0}^d m_{k}P_{ki}.$$

Using equality (\ref{crossPQ}), we obtain
$$\langle F,E_j\rangle =
\frac{1}{v}\sum_{i=0}^d Q_{ij}\sum_{k=0}^d P_{ki}m_k =
\sum_{k=0}^d m_k \sum_{i=0}^d \frac{Q_{ij}P_{ki}}{v} =
\sum_{k=0}^d m_k \delta_{j,k} = m_j,$$ which proves the theorem.
\medskip

Therefore it follows from Theorem 2 that, for a putative equitable
partition of association scheme, condition (\ref{conditionGodsil})
cannot say more than Theorem 1. Moreover, the following example
shows that the condition may be even weaker than an obvious
condition of integrality of the elements of quotient matrix.

Consider a perfect matching $M=\{(i,i'): i \in \{0,\ldots,4\}\}$.
Let $C$ be a labeled cycle with vertex set $\{0,\ldots,4\}$, and
$C'$ be its complement with vertex set $\{0',\ldots,4'\}$ so that
$i'\sim_{C'} j'$ if and only if $i\not\sim_C j$.

Define a graph $\Gamma$ with vertex set
$\{i:i \in \{0,\ldots,4\}\} \cup \{i': i \in \{0,\ldots,4\}\}$
and edge set consisting of edges $C$, $C'$, and $M$. It is easy to see
that $\Gamma$ is the Petersen graph, which is known to be
{\it distance-regular} with diameter 2, i.e., the distance relations on
its vertex set form an association scheme with two classes.

Let us consider a 5-partition consisting of a cell with a pair of
adjacent vertices $\Gamma$ (say, $C_1=\{0,0'\}$) and of four cells of
pairs of non-adjacent vertices (say, $C_2=\{1,2'\}$,
$C_3=\{2,1'\}$, $C_4=\{3,4'\}$, $C_5=\{4,3'\}$). Clearly, the
partition is not equitable, because, for example, the vertex $2'$
has a neighbour in $C_4$, whereas the vertex 1 does not.

However, the partition is feasible with respect to (\ref{conditionGodsil}).
It is easy to see that $\langle F,A_0\rangle$ is always the number of cells
in the partition, i.e. $$\langle F,A_0\rangle=5,$$ and, further,
$$\langle F,A_1\rangle=1, \langle F,A_2\rangle=4.$$

Using equality (\ref{conditionGodsil}), we obtain:
$$\langle F,E_0\rangle=1,  \langle F,E_1\rangle=2, \langle F,E_2\rangle=2,$$
which completes our example.


\bigskip
\begin{thebibliography}{1}

\bibitem{Ban84} E. Bannai, T. Ito, {\it Algebraic combinatorics. I},
Benjamin/Cummings Publishing Co. Inc., Menlo Park, CA, 1984.

\bibitem{Bra89}
A.E. Brouwer, A.M. Cohen, A. Neumaier, {\it Distance regular
graphs}, Berlin: Springer-Verl., 1989.

\bibitem{Cam} P.J. Cameron, {\it Permutation groups}, London Math. Soc.
Student Texts N45. Cambridge: Cambridge Univ. Press. 1999.

\bibitem{Del73}  P. Delsarte, {\it An algebraic approach to the association schemes of coding theory}, Philips Res. Rep. Suppl., {\bf 10} (1973), 1--97.

\bibitem{God04}
C. Godsil, R. Gordon, {\it Algebraic graph theory}, Springer Science+Business 
Media, LLC, 2004.

\bibitem{God10} C. Godsil, {\it Association schemes}, University of Waterloo, 2010.

\bibitem{Makh10} A.A. Makhnev, {\it On automorphisms of distance-regular graphs}, 
J. Math. Sci. (New York), {\bf 166}:6 (2010),  733--742.

\bibitem{Martphd} W.J. Martin, {\it Completely regular subsets}, Ph.D. thesis, 
University of Waterloo, 1992.

\bibitem{Van11}
F. Vanhove, {\it Incidence geometry from an algebraic graph theory
point of view}, PhD Thesis, University of Ghent, 2011.

\end{thebibliography}
\end{document}